\documentclass[10pt,a4paper,oneside,bibnotes]{article}

\usepackage{amsfonts}
\usepackage[pdftex,dvipsnames]{xcolor}
\usepackage{url}
\usepackage{amsmath}
\usepackage{hyperref}
\usepackage{amsthm}
\usepackage{graphicx}
\usepackage{amscd}
\usepackage{amssymb}
\usepackage{latexsym}
\usepackage[utf8]{inputenc}
\usepackage{theoremref}
\usepackage[T1]{fontenc}
\usepackage{lmodern}
\usepackage{verbatim}
\usepackage{float}
\usepackage[a4paper]{geometry}
\usepackage{parskip}
\usepackage{mathrsfs}
\usepackage{tikz}
\usepackage{lipsum}        
\usepackage{xargs}           
\usepackage[colorinlistoftodos,prependcaption,textsize=tiny]{todonotes}
\usepackage{tikz-cd}
\usetikzlibrary{matrix}
\usepackage{setspace}
\usepackage{pgfkeys}              
\usepackage{xfrac}
\usepackage{mathtools}
\usepackage{titlesec}
\usetikzlibrary{fit, arrows, calc, positioning}
\theoremstyle{plain}
\usepackage[utf8]{inputenc}
\newcommand*\rfrac[2]{{}^{#1}\!/_{#2}}
\newcommand{\Q}{{\mathbb Q}}
\newcommand{\Z}{{\mathbb Z}}
\newcommand{\F}{{\mathbb F}}

\newcommand{\C}{{\mathbb C}}
\newcommand{\ra}{{\rightarrow}}

\newcommand{\Gn}{G_{n}}
\newcommand{\Go}{G_{0}}

\newcommand{\Hn}{\Delta_{n}}
\newcommand{\Bo}{B_{0}}
\newcommand{\Bn}{B_{n}}
\newcommand{\Po}{{\PR^{1}_{0}}}
\newcommand{\Pn}{{\PR^{1}_{n}}}

\newcommand{\PR}{{\mathbb P}}

    \DeclareFontFamily{U}{wncy}{}
    \DeclareFontShape{U}{wncy}{m}{n}{<->wncyr10}{}
    \DeclareSymbolFont{mcy}{U}{wncy}{m}{n}
    \DeclareMathSymbol{\Sh}{\mathord}{mcy}{"58} 
\DeclareMathOperator{\Gal}{Gal}

\DeclareMathOperator{\Sel}{Sel}
\DeclareMathOperator{\PGL}{PGL}
\DeclareMathOperator{\Fr}{Fr}
\DeclareMathOperator{\Sym}{Sym}
\DeclareMathOperator{\PSL}{PSL}
\DeclareMathOperator{\im}{im}
\DeclareMathOperator{\Ind}{Ind}
\DeclareMathOperator{\Hom}{Hom}
\DeclareMathOperator{\Nm}{Nm}
\DeclareMathOperator{\HC}{H}
\DeclareMathOperator{\ord}{ord}
\DeclareMathOperator{\SL}{SL}
\DeclareMathOperator{\GL}{GL}

\DeclareMathOperator{\End}{End}
\DeclareMathOperator{\Irr}{Irr}
\DeclareMathOperator{\Aut}{Aut}

\DeclareMathOperator{\Res}{Res}
\DeclareMathOperator{\disc}{disc}

\DeclareMathOperator{\Ext}{Ext}

\def\m{{\ensuremath{\mathfrak{m}}}}

\DeclareMathOperator{\M}{M}

\DeclareMathOperator{\tr}{tr}
\DeclareUnicodeCharacter{2212}{-}
\makeatletter
\def\thm@space@setup{%
  \thm@preskip=\parskip \thm@postskip=0pt
}
\titleformat{\subsection}
  {\normalfont\bfseries\fontsize{8}{8}\centering}{\thesubsection}{1em}{}

\titleformat{\section}
  {\normalfont\bfseries\fontsize{14}{15}\centering}{\thesection}{1em}{}
\makeatother

\title{Kolyvagin Derivatives of Modular Points on Elliptic Curves}
\author{Richard Hatton}
\date{}

\begin{document}

\newtheorem*{theorem*}{Theorem}
\newtheorem{theorem}{Theorem}[section]
\newtheorem{lemma}[theorem]{Lemma}
\newtheorem{definition}[theorem]{Definition}
\newtheorem{proposition}[theorem]{Proposition}
\newtheorem{corollary}[theorem]{Corollary}
\newtheorem{conjecture}[theorem]{Conjecture}
\newtheorem{remark}{Remark}
\newtheorem*{remark*}{Remark}
\newtheorem{example}{Example}
\renewcommand\theexample{\unskip}

\maketitle

\begin{abstract}
Let $E/\Q$ and $A/\Q$ be elliptic curves. We can construct modular points derived from $A$ via the modular parametrisation of $E$. With certain assumptions we can show that these points are of infinite order and are not divisible by a prime $p$. In particular, using Kolyvagin’s construction of derivative classes, we can find elements in certain Shafarevich-Tate groups of order $p^{n}$.
\end{abstract}

\bigskip

\section*{INTRODUCTION}

Let $E/\Q$ be an elliptic curve of conductor $N_{E}$. Then due to the modularity theorem, there exists a surjective morphism  \[\phi_{E}:X_{0}(N_{E}) \ra E\] defined over $\Q$ known as the modular parametrisation of $E$, where $\infty$ on the modular curve $X_{0}(N_{E})$ is mapped to $O$. There exists a subvariety $Y_{0}(N_{E})$ of $X_{0}(N_{E})$ which is a moduli space of points $x_{A,C}=(A,C)$ where $A$ is an elliptic curve and $C$ is a cyclic subgroup of $A$ of order $N_{E}$. Fixing $A/\Q,$ the image of $x_{A,C}$ under $\phi_{E}$ is known as a \emph{modular point}, which we will denote $P_{A,C} \in E(\Q(C))$ where $\Q(C)$ is the field of definition of $C$. We denote the compositium of all such $\Q(C)$ as $K_{N_{E}}$. This is the smallest field $K$ such that its absolute Galois group $G_{K}$ acts by scalars on $A[N_{E}].$ This has Galois group $G_{N_{E}}:=\Gal(K_{N_{E}}/\Q)$ which can be identified as a subgroup of $\PGL_{2}(\Z/N_{E}\Z)$ via the mapping \[\overline{\tau}_{\Q,A,N_{E}}:G_{\Q} \ra \Aut(A[N_{E}])\cong \GL_{2}(\Z/N_{E}\Z) \ra \PGL_{2}(\Z/N_{E}\Z).\]

We can also define \emph{higher modular points} above $P_{A,C}$. These are points of the form $\phi_{E}(B,D)$ for an elliptic curve $B$ isogenous to $A$ over $\overline{\Q}$ and $D$ a cyclic subgroup of $B$ of order $N_{E}$.

In this paper, we use these points to bound Selmer groups using methods similar to that of Kolyvagin in \cite{VK90} and Wuthrich in \cite{CW09}. Kolyvagin initially looked at a specific type of modular point known as Heegner points and used them to bound Selmer groups. This involved creating cohomology classes coming from these points and using the classes to bound the Selmer groups from above. Wuthrich then worked on an analogue system to Kolyvagin's work where he uses a type of modular point known as self points to create derivative classes and finds lower bounds of Selmer groups over certain fields.

We look to extend the idea of self points to a generalised construction of modular point. In particular, we use more advanced methods in modular representation theory to show that Selmer groups over certain fields must contain points of prime power order when the higher modular points satisfy certain conditions.

In the first section, we look at the divisibility of the modular points in $E(\Q(C))$ and understand the relationships between them. This will provide us with an understanding of the rank of the group generated by these points. Initially, we want to see when the modular points are of infinite order. We obtain the following result.

\begin{theorem*}
Let $E/\Q$ be an elliptic curve of conductor $N_{E}$. Let $A/\Q$ be an elliptic curve such that the $j$ invariant of $A$ is not in $\frac{1}{2}\Z$ and the degree of any isogeny of $A$ defined over $\Q$ is coprime to $N_{E}$. Then the modular points $P_{A,C}$ are of infinite order for all cyclic subgroups $C$ of order $N_{E}$ in $A$.
\end{theorem*}

From this, we can show that if $p$ is a prime with specific conditions related to $A$ and $E$, then $P_{A,C}$ is not divisible by $p$ in $E(\Q(C))$ as seen in \thref{p div}. We also see that there exist relationships between the modular points. If $d$ is a divisor of $N_{E}$ and $B$ is a cyclic subgroup of $A$ of order $d$, then \[\sum\limits_{C \supset B}P_{A,C}\in E(K_{d})\text{ is torsion.}\] This reduces the rank of the group generated by these points. 

We then take a look at a specific case of creating higher modular points for a prime $p$ of either good ordinary or multiplicative reduction with respect to $E$. Here, we will look at the case where $p$ is coprime to $N_{E}$. Let $D$ be a cyclic subgroup of $A$ of order $p^{n+1}$ for a prime $p$ and $n \ge 0$. We look at the higher modular point coming from $(A/D,\psi(C))$ where $\psi:A \ra A/D$ is the isogeny defined by $D$. We define $Q_{A,D}=\phi_{E}(A/D,\psi(C)) \in E(\Q(C,D))$. We see that the higher modular points form a trace-compatible system with \[a_{p}(E) \cdot Q_{A,D}=\sum\limits_{D' \supset D}Q_{A,D'}\] where the sum is taken over the subgroups $D'$ of $A$ of order $p^{n+2}$ containing $D$ and $a_{p}(E)$ is the $p$-th Fourier coefficient for the modular form associated to the isogeny class of $E$. Using this, we can show that if $P_{A,C}$ is of infinite order, then we have created a tower of points which are also of infinite order. In particular, we show that the higher modular points generate a group of rank $p^{n+1}+p^{n}$ as seen in \thref{good} and so if $F_{n}$ defines the compositium of all such $\Q(C,D)$, then rank($E(F_{n})) \ge p^{n+1}+p^{n}$ by \cite[Corollary 3.7]{MH79}. Hence we are able to use this information to establish a link between the group generated by the higher modular points and the representation theory associated to the projective general linear group. 

If we let $\Gn:=\PGL_{2}(\Z/p^{n+1}\Z)$ for $n \ge 0$, we will see that we can relate the group generated by the higher modular points to certain $\Z_{p}[\Gn]$-lattices in 
\[V_{n}=\ker\Big(\Q_{p}\big[\rfrac{\Gn}{\Bn}\big] \ra \Q_{p}\Big)\] for $p$ prime and $\Bn$ a Borel subgroup of $\Gn$. We can view \[V_{n}=\bigg\{f:\Pn \ra \Q_{p}:\sum\limits_{C}f(C)=0\bigg\}\] where $\Pn:=\PR^{1}(\Z/p^{n+1}\Z)$. This contains the standard lattice \[T_{n}=\ker\Big(\Z_{p}\big[\rfrac{\Gn}{\Bn}\big] \ra \Z_{p}\Big),\] which we can easily understand the cohomology of with respect to subgroups of $\Gn$. Hence, we take a look at the integral representation theory of $\Gn$ and will later look at the application of this to the modular points. 

We then follow a similar method to Wuthrich in \cite{CW09} from the ideas of Kolyvagin in \cite{VK90}. We create derivative classes coming from higher modular points of infinite order which are not divisible by $p$ in $E(\Q(C))$ as shown in \thref{p div}.

Suppose $p$ is one of the following:
	\begin{itemize}
	\item A prime of split multiplicative reduction for $E$ and $p \nmid \ord_{p}(\Delta_{E})$;
	\item A prime of non-split multiplicative reduction for $E$;
	\item A non-anomalous prime of good ordinary reduction for $E$;
	\end{itemize}
where $\Delta_{E}$ is the minimal discriminant of $E$. This ensures that the higher modular points are not divisible by $p$ in $E(\Q(C))$. Let 
\[
    \ F_{n}:= \left\{
                \begin{array}{ll}
                  K_{p^{n+1}N_{E}} & \text{if $p\nmid N_{E}$},\\
                  K_{p^{n}N_{E}} & \text{if $p || N_{E}$},
                \end{array}
              \right.
\] for $n \ge -1$ with $F:=F_{-1}$. We assume \[\tau_{F,A,p}:\Gal(\overline{F}/F) \ra \PGL_{2}(\Z_{p})\] is surjective giving $\Gal(F_{n}/F) \cong \Gn$. We let $A_{n}$ be a non-split Cartan subgroup of $\Gn$. This is a cyclic subgroup of order $p^{n+1}+p^{n}$. Then we define $L_{n}$ to be the subfield of $F_{n}$ fixed by $A_{n}$. We are able to construct a mapping \[\delta_{n}:\HC^{1}(A_{n},S_{n}) \ra \Sh(E/L_{n})[p^{\infty}]\] where $S_{n}$ denotes the saturated group generated by the higher modular points in $E(F_{n})$. With the conditions on $p$, we can show that the source of $\delta_{n}$ is a cyclic group of order $p^{n}$. This result derives from the link between the construction of the derivative classes and integral representation theory. The group $S_{n}$ defined earlier is isomorphic to a $\Z_{p}[\Gn]$-lattice containing $T_{n}$. Due to the structure of $S_{n}$, we are able to show that $S_{n} \cong T_{n}$ under the conditions we have stated and as we understand the cohomology of $T_{n}$ associated to subgroups of $\Gn$. This leads to the following.

\begin{theorem*}
Let $E/\Q$ and $A/\Q$ be elliptic curves of conductor $N_{E}$ and $N_{A}$ respectively. Let $p>3$ be a prime of multiplicative reduction for $A$ such that $\ord_{p}(\Delta_{A})=1$ and let $F_{n}$ be as defined above. Assume that:
	\begin{enumerate}
	\item $A$ is semistable;
	\item $E$ has either split multiplicative reduction at $p$ with $p \nmid \ord_{p}(\Delta_{E})$, non-split multiplicative reduction at $p,$ or good ordinary non-anomalous reduction at $p$;
	\item The degree of any isogeny of $A$ defined over $\Q$ is coprime to $N_{E}$;
	\item $\overline{\rho}_{\Q,A,p}:G_{\Q} \ra \Aut(A[p]) \cong \GL_{2}(\F_{p})$ is surjective;
	\item $\overline{\rho}_{\Q,E,p}:G_{\Q} \ra \Aut(E[p]) \cong \GL_{2}(\F_{p})$ is surjective;
	\item Any prime $\ell$ of bad reduction for $E$ and good reduction for $A$ has $a_{\ell}(A)^{2}-4\ell$ square modulo $p$.
	\end{enumerate}
Then there exists an element of order $p^{n}$ in $\Sel^{p^{n}}(E/L_{n})$.
\end{theorem*}

We are able to prove that if the first and fourth conditions in the theorem are true, then $\tau_{F,A,p}$ is surjective and so $\Gal(F_{n}/F) \cong G_{n}$ for all $n \ge 0$. The fifth point is essential in order to ensure $E(F_{n})[p]=0$ for all $n \ge 0$. We know that this condition, together with the fourth condition, excludes only a finite number of primes if $A$ and $E$ have no complex multiplication. Therefore, these two conditions aren't very strong.  The third condition is vital to ensure all the points $P_{A,C}$ are of infinite order. This condition is weak as the number of such isogenies is finite. 

The final point ensures that the primes dividing $N_{E}$ but not $N_{A}$ split completely in $F_{n}/L_{n}$. This is fundamental when showing that the image of the derivation map $\delta_{n}$ lies in $\Sh(E/L_{n})[p^{\infty}]$. In particular, as this is only focusing on primes of bad reduction for $E$ and good reduction for $A$, this condition only looks at a finite number of primes. 

However, we are unsure as to whether the image of this element in $\Sel^{p^{n}}(E/L_{n})$ is trivial or not in $\Sh(E/L_{n})[p^{n}]$. This will require a more in-depth look at the derivation map $\delta_{n}$ to see whether this mapping is injective or not. Also, we still do not fully understand all potential $\Z_{p}[\Gn]$-lattices that $S_{n}$ could be isomorphic to. Further research into the integral representation theory of $\Z_{p}[\Gn]$ would improve our understanding of the structure of the saturated group of higher modular points and further still, understand the properties of the derivative classes constructed.

\section{PRELIMINARIES} 

Let $K$ be a number field. For an elliptic curve $E$ over $K$ and $m>1$ an integer, we let $E[m]$ be the $m$-torsion subgroup of $E(\overline{K})$. We have $G_{K}$ acting on $E[m]$ where $G_{K}:=\Gal(\overline{K}/K)$ is the absolute Galois group of $K$. This leads to a Galois representation \[\overline{\rho}_{K,E,m}:G_{K} \ra \Aut(E[m]) \cong \GL_{2}(\Z/m\Z).\] 

Let $T_{p}E=\varprojlim\limits_{n}E[p^{n}]$ be the $p$-adic Tate module of $E$ for a prime $p$. Then $G_{K}$ acts on $T_{p}E$ which leads to the Galois representation \[\rho_{K,E,p}:G_{K} \ra \Aut(T_{p}E) \cong \GL_{2}(\Z_{p}).\]  We define the mapping $\overline{\tau}_{K,E,m}:= s_{m} \circ \overline{\rho}_{K,E,m}$ where $s_{m}$ is the quotient mapping to $\PGL_{2}(\Z/m\Z)$ and the mapping $\tau_{K,E,p}:= s \circ \rho_{K,E,p}$ where $p$ is prime and $s$ is the quotient mapping to $\PGL(\Z_{p})$. Throughout, we will denote the centre of $\GL_{2}(\Z/m\Z)$ as $Z_{m}$ and define $\Gn:=\PGL_{2}(\Z/p^{n+1}\Z)$ for a prime $p$ and $n \ge 0$. 

In this article, we will be looking at the links between the modular points we have constructed and the representation theory associated to $\Gn$. We let $\Bn$ denote a Borel subgroup of $\Gn$. We see that $\Gn$ acts on the projective line over $\Z/p^{n+1}\Z$ for $n \ge 0$ via linear substitution, which will be denoted $\Pn$ throughout.

\section{MODULAR POINTS ON ELLIPTIC CURVES}

In this section, we will see when the points $P_{A,C}$ are of infinite order. Further, when the prime $p$ satisfies certain conditions with respect to $A$ and $E$, we will show that the points $P_{A,C}$ are not divisible by $p$ in $E(\Q(C))$. We prove the following.

\begin{theorem}
Let $A/\Q$ and $E/\Q$ be elliptic curves of conductor $N_{A}$ and $N_{E}$ respectively. Suppose the $j$ invariant of $A$ is not in $\frac{1}{2}\Z$. Then there exists a cyclic subgroup of order $N_{E}$ in $A$, denoted $C$, such that $P_{A,C} \in E(\overline{\Q}_{p})$ is non-torsion.
\end{theorem}

\begin{proof}
The following has been adapted from the proof of \cite[Theorem 2]{CW09} for the more general setting. Let $p$ be a prime which divides the denominator of the $j$-invariant of $A$. If $p^{2} \mid N_{A}$, we know $A$ acquires multiplicative reduction at $p$ over some extension of $\Q$. Fix an embedding of $\overline{\Q}$ into $\overline{\Q}_{p}$. We consider the modular parametrisation over $\overline{\Z}_{p}$. The modular curve $X_{0}(N_{E})$ over $\overline{\Z}_{p}$ has a neighbourhood of the cusp $\infty$ consisting of couples $(J,C)$ of the Tate curve of the form $J(\overline{\Q}_{p})=\overline{\Q}_{p}^{\times}/q^{\Z}$ together with a cyclic subgroup $C$ of order $N_{E}$ generated by the $N_{E}^{\text{th}}$ root of unity. The parameter $q$ is a $p$-adic analytic uniformiser at $\infty$, so that the Spf $\overline{\Z}_{p}[[q]]$ is the formal completion of $X_{0}(N_{E})/\overline{\Z}_{p}$ at the cusp $\infty$, as seen in \cite[Chapter 8]{KM85}.

Since $A$ has multiplicative reduction over $\overline{\Z}_{p}$, there is exactly one point $x_{A,C}$ in the neighbourhood of $\infty$ on $X_{0}(N_{E})$ which is represented by the $p$-adic Tate parameter $q_{A}$ associated to $A$ and $C$ is isomorphic to $\mu[N_{E}]$. Let $f_{E}=\sum a_{n}(E)q^{n}$ be the normalised newform associated to $E$ and let $\omega_{E}$ be the invariant differential on $E$. Then we have $\phi_{E}^{*}(\omega_{E})=c_{E} \cdot f_{E}/q \cdot dq$ where $c_{E}$ is the Manin constant of $E$. Hence \[\log_{E}(\phi_{E}(x_{A,C}))=\int^{\phi_{E}(x_{A,C})}_{O}\omega_{E}=c_{E} \cdot \int^{q_{A}}_{0}f_{E}\frac{dq}{q}=c_{E} \cdot \sum\limits_{n=1}^{\infty}\frac{a_{n}(E)}{n} \cdot q^{n}_{A}\] where $\log_{E}$ is the formal logarithm associated to $E$ from the formal group $\hat{E}(\m)$ to the maximal ideal $\m$ of $\overline{\Z}_{p}$. We let $|\cdot |_{p}$ be the normalised absolute value such that $|p|_{p}=p^{-1}$. We then have the following lemma.

\begin{lemma} \cite[Lemma 3]{CW09}
Let $(J,C)$ be a point on $Y_{0}(N_{E})(\overline{\Q}_{p})$ such that $J$ is isomorphic to the Tate curve with parameter $q_{0} \neq 0$ and $C$ is isomorphic to $\mu[N_{E}]$. If $|q_{0}|_{p}<p^{-\frac{1}{p-1}}$ then $\phi_{E}(J,C)$ is a point of infinite order on $E(\overline{\Q}_{p})$.
\end{lemma}

Hence, we have for $p \neq 2$ that \[|q_{A}|_{p}=|j(A)|_{p}^{-1} \leq p^{-1}<p^{-\frac{1}{p-1}}\] and if $p=2$ then \[|q_{A}|_{2}=|j(A)|_{2}^{-1} \leq p^{-2}<p^{-\frac{1}{p-1}}.\] Therefore, by the lemma, we have that $P_{A,C}=\phi_{E}(x_{A,C})$ is a point of infinite order in $E(\overline{\Q}_{p})$.
\end{proof}

We assume for primes $\ell$ of bad reduction for $E$, there exists no $\ell$-isogeny on $A$. Hence the set of modular points $\{P_{A,C}\}_{C}$ form a single orbit under the action of $G_{N_{E}}.$ This shows that the points $P_{A,C}$ are not equal and have the same order for all $C$. Hence, $P_{A,C}$ is of infinite order in $E(K_{N_{E}})$ for all $C$.

We can also show that there exist relationships between these points. These are in fact identical to the relations for the self point case as in \cite[Proposition 4]{CW09}.

\begin{proposition}
The sum of the modular points $P_{A,C}$ as $C$ varies through all cyclic subgroups of $A$ of order $N_{E}$ is a torsion point defined over $\Q$. Let $d \neq N_{E}$ be an integer dividing $N_{E}$. Then \[\sum\limits_{C \supset B}P_{A,C} \in E(K_{d})\] is a torsion point where $B$ is a cyclic subgroup of $A$ of order $d$ and the sum is taken over all cyclic subgroups of $A$ of order $N_{E}$ containing $B$.
\end{proposition}

\begin{proof}
Identical to the proof of \cite[Proposition 4]{CW09}.
\end{proof}

Hence, we have provided the terms in which the points are of infinite order and have found the relationships between them in order to understand the group generated by these points. We can in fact further show that the points are not divisible by $p$ in $E(\Q(C))$ when the prime $p$ satisfies certain conditions.

\begin{proposition} \thlabel{p div}
Suppose $j(A) \notin \frac{1}{2}\Z$ and $p>2$ is a prime of multiplicative reduction for $A$ such that it is one of the following:
	\begin{itemize}
	\item A prime of split multiplicative reduction for $E$ and $p \nmid \ord_{p}(\Delta_{E})$;
	\item A prime of non-split multiplicative reduction for $E$;
	\item A non-anomalous prime of good reduction for $E$.
	\end{itemize}
Let $\Delta_{A}$ be the minimal discriminant for $A$. Then $P_{A,C} \notin p^{\ord_{p}(\Delta_{A})} \cdot E(\Q(C))$.
\end{proposition}

\begin{proof}
Let $\mathfrak{p}$ be the place of $\Q(C)$ corresponding to the chosen embedding $\overline{\Q}$ to $\overline{\Q}_{p}$. Then we have $\Q(C)_{\mathfrak{p}}=\Q_{p}$. As $p$ is one of the primes in the proposition with respect to $E$, then the order of the group of components of $E$ over $\Q_{p}$ and the number of non-singular points in the reduction $\tilde{E}(\F_{p})$ are both prime to $p$. Hence, if $P_{A,C}$ is divisible by $p^{\ord_{p}(\Delta_{A})}$ in $E(\Q_{p})$ then it is divisible by $p^{\ord_{p}(\Delta_{A})}$ in $\hat{E}(p\Z_{p})$. But the valuation of $\log_{E}(P_{A,C})$ shows this cannot happen. Let $z=\log_{E}(P_{A,C})$. As $\Q(C)_{\mathfrak{p}}=\Q_{p}$, then we have $\ord_{p}(q_{A})=\ord_{p}(\Delta_{A})$. We then see that
\begin{align*}
\ord_{p}(z)	&=\min\{\ord_{p}(q_{A}),\ord_{p}(a_{p}(E)) +
			p \cdot \ord_{p}(q_{A})-\ord_{p}(p)\} \\
			&=\min\{\ord_{p}(\Delta_{A}),\ord_{p}(a_{p}(E)) +p \cdot \ord_{p}(\Delta_{A})-1\}=\ord_{p}(\Delta_{A})
\end{align*}
as $a_{p}(E)$ is an integer. Therefore $P_{A,C} \in \hat{E}(p^{\ord_{p}(\Delta_{A})}\Z_{p})$. If $p \neq 2$, then for all $n \ge 2$, we have $n-\ord_{p}(n) \ge 2$. Therefore, we have
\begin{align*}
\ord_{p}(z-q_{A})  	& \ge \min\limits_{n \ge 2}\big[\ord_{p}(a_{n}(E))-\ord_{p}(n)+n \cdot \ord_{p}(q_{A})\big]\\
			&\ge \min\limits_{n \ge 2}\big[n \cdot \ord_{p}(\Delta_{A})-\ord_{p}(n)\big]\\
			&\ge \min\limits_{n \ge 2}\big[(n-\ord_{p}(n)) \cdot \ord_{p}(\Delta_{A})\big]\\
			&\ge 2\ord_{p}(\Delta_{A}).
\end{align*}
Hence, $z$ is congruent to $q_{A}$ modulo $p^{2\ord_{p}(\Delta_{A})}$. So the value $z$ is non zero and therefore ${P_{A,C} \notin p^{\ord_{p}(\Delta_{A})} \cdot E(\Q(C))}$.
\end{proof}

\section{HIGHER MODULAR POINTS} \label{s3}

We now extend the construction of modular points to higher modular points by following the construction as in \cite{CW09} in the more generalised setting. We will create a tower of number fields in which the higher modular points will be defined and will show that they satisfy certain trace relations similar to that in Kolyvagin's construction of Heegner points. In particular, the relations between the points show that a link can be established with representations of $\Gn$.

Let $E/\Q$ and $A/\Q$ be elliptic curves of conductor $N_{E}$ and $N_{A}$ respectively and let $D$ be a cyclic subgroup in $A$. We then take the isogenous curve $A/D$ and a cyclic subgroup of order $N_{E}$ in this isogenous curve and look at the image of the point under $\phi_{E}$. We will focus on two cases of construction of higher modular points. The first will be the multiplicative reduction case and the second the good reduction case.

{\footnotesize \subsection{MULTIPLICATIVE CASE}}

Let $E/\Q$ be an elliptic curve of conductor $N_{E}$ and let $p$ be a prime of multiplicative reduction. We let $N_{E}=pM$ and define $F_{n}:=K_{p^{n+1}M}$ for all $n \ge -1$ with $F:=F_{-1}$. We can view $\Hn:=\Gal(F_{n}/F)$ as a subgroup of $\Gn$.

Let $B$ be a cyclic subgroup in $A$ of order $M$. Let $n \ge 0$ and $D$ be a cyclic subgroup of order $p^{n+1}$ in $A$. Let $C=D[p] \oplus B$ which is a cyclic subgroup of order $N_{E}$ and let $\psi:A \ra A/D$ be the isogeny associated with $D$ with dual $\hat{\psi}$. We define \[C'=\ker(\hat{\psi})[p] \oplus \psi(B).\] This is a cyclic subgroup in $A/D$ of order $N_{E}$. We then define $y_{A,D}=(A/D,C') \in Y_{0}(N_{E})$. We define a higher modular point above $P_{A,C}$ as the point $Q_{A,D}=\phi_{E}(y_{A,D}) \in E(F_{n}).$

When $n=0$, we see from above that we have $y_{A,D}=w_{p}(x_{A,C})$ where $w_{p}$ is the Atkin-Lehner involution on $X_{0}(N_{E})$. Therefore, we have $Q_{A,D}=-a_{p}(E) \cdot P_{A,C} +T$ for some $T \in E(\Q)[2]$ and $a_{p}(E) \in \{\pm 1\}$ is the Hecke eigenvalue of the newform $f$ attached to $E$.

Let $D$ be a cyclic subgroup of order $p^{n+1}$ in $A$. Let $T_{p}$ be the Hecke operator $T_{p}$ on $J_{0}(N_{E})$. Then \[T_{p}((y_{A,D})-(\infty))=\sum\limits_{D' \supset D}((y_{A,D'})-(\infty))\] where the sum runs over all cyclic subgroups $D'$ of order $p^{n+2}$ in $A$ containing $D$. This then leads to the trace relation \[a_{p}(E) \cdot Q_{A,D}=\sum\limits_{D'\supset D}Q_{A,D'}.\] We know that if $J,J' \supset D$ where $J$ and $J'$ are cyclic subgroups of order $p^{n+2}$, then the points $Q_{A,J}$ and $Q_{A,J'}$ lie in the same Galois orbit with trivial stabiliser. In particular, this ensures that these points are not equal and have the same order. Hence, by induction, if $P_{A,C}$ is of infinite order then so is $Q_{A,D}$. Also, if $P_{A,C} \notin p \cdot E(F_{0})$ then $Q_{A,D} \notin p \cdot E(F_{n})$. 

We are able to obtain the following lemma, which is a generalised version of \cite[Lemma 16]{CW09} and is proven in an identical way.

\begin{lemma}
Fix a cyclic subgroup $B$ in $A$ of order $M$ and let $n \ge 0$. Then $\sum\limits_{D}Q_{A,D}$ is a torsion point in $E(F)$ where the sum is over all cyclic subgroups $D$ in $A$ of order $p^{n+1}$.
\end{lemma}

Let \[V_{n}=\ker\big(\Q_{p}\big[\rfrac{\Gn}{\Bn}\big]\ra \Q_{p}\big)\] be a $\Q_{p}[\Hn]$-module of dimension $p^{n+1}+p^{n}-1$. We can view \[V_{n}=\bigg\{f:\Pn \ra \Q_{p}:\sum\limits_{C}f(C)=0\bigg\}\] where $\Gn$ acts on $\Pn$ by permutating the basis. Then if we fix a cyclic subgroup $B$ of order $M$ in $A$, we have a $\Hn$-morphism \[V_{n} \ra E(F_{n}) \otimes \Q_{p}\]\[e_{D} \mapsto Q_{A,D}\] where $e_{D}$ is the characteristic function on $D$. If we then assume that $\tau_{F,A,p}$ is surjective, then we have $\Hn \cong \Gn$. We thus have the following lemma, a generalised version of \cite[Theorem 17]{CW09} and is proven in an identical way.

\begin{proposition}
Let $E/\Q$ and $A/\Q$ be elliptic curves of conductor $N_{E}$ and $N_{A}$ respectively. Let $p$ be a prime of multiplicative reduction of $E$ and let $P_{A,C} \in E(F_{0})$ have infinite order. Suppose $\tau_{F,A,p}$ is surjective. Then all higher modular points $Q_{A,D}$ above $P_{A,C}$ are of infinite order and generate a group of $\Q_{p}$-dimension $p^{n+1}+p^{n}-1$ in $E(F_{n}) \otimes \Q_{p}$.
\end{proposition}

{\footnotesize \subsection{GOOD CASE}}

Let $p$ be a prime of good reduction for $E$. We define $F_{n}:=K_{p^{n+1}N_{E}}$ where $F:=F_{-1}$ is a number field such that $P_{A,C}$ is of infinite order in $E(F)$. Hence, we can view $\Hn=\Gal(F_{n}/F)$ as a subgroup of $\Gn$. 

For $n \ge 0$, we let $D$ be a cyclic subgroup in $A$ of order $p^{n+1}$ and then we can construct the higher modular points. We let $\psi:A \ra A/D$ be the isogeny associated to $D$. Then we have $y_{A,D}=(A/D,\psi(C)) \in Y_{0}(N_{E}).$ We then define the higher modular point above $P_{A,C}$ as \newline $Q_{A,D}=\phi_{E}(y_{A,D}) \in E(F_{n})$. 

As before, if we use the Hecke operator $T_{p}$, we obtain the trace relation \[a_{p}(E) \cdot Q_{A,D}=\sum\limits_{D' \supset D}Q_{A,D'}\] for all $n \ge 0$ and $D$ as defined before. Here, the sum runs over all cyclic subgroups $D'$ of order $p^{n+2}$ containing $D$. Hence, we have \[a_{p}(E) \cdot P_{A,C}=\sum\limits_{D}Q_{A,D}\] where the sum runs over all cyclic subgroups of order $p$. Therefore, as in the multiplicative case, if $J,J' \supset D$ where $J$ and $J'$ are cyclic subgroups of order $p^{n+2}$, then the points $Q_{A,J}$ and $Q_{A,J'}$ are not equal and have the same order. Therefore, if $P_{A,C}$ is of infinite order in $E(F)$, then $Q_{A,D}$ is of infinite order in $E(F_{n})$.

Therefore, if we let $V_{(n)}=\Q_{p}\big[\rfrac{\Gn}{\Bn}\big]=$ Maps$(\Pn \ra \Q_{p})$ be a $\Q_{p}[\Hn]$-module of dimension $p^{n+1}+p^{n}$, then we can define a $\Hn$-morphism \[V_{(n)} \ra E(F_{n}) \otimes \Q_{p}\]\[e_{D} \mapsto Q_{A,D}.\] If we then assume that $\tau_{F,A,p}$ is surjective, then we have $\Hn \cong \Gn$. Again, we are able to obtain the following lemma, which is a generalised version of \cite[Theorem 18]{CW09} and is proven in an identical way.

\begin{proposition} \thlabel{good}
Let $E/\Q$ and $A/\Q$ be elliptic curves of conductor $N_{E}$ and $N_{A}$ respectively. Let $p$ be a prime of good ordinary reduction of $E$ ensuring $a_{p}(E) \neq 0$ and let $P_{A,C} \in E(F)$ have infinite order. Suppose $\tau_{F,A,p}$ is surjective. Then all higher modular points $Q_{A,D}$ above $P_{A,C}$ are of infinite order and generate a group of $\Q_{p}$-dimension $p^{n+1}+p^{n}$ in $E(F_{n}) \otimes \Q_{p}$.
\end{proposition}

We can thus understand the group generated by the higher modular points by studying the representations of $\Gn$ associated to them.

\section{REPRESENTATIONS OF {\normalsize $\PGL_{2}(\Z/p^{n+1}\Z)$}}

We are interested in the representations that appear in the study of the higher modular points. We have seen that when $\tau_{F,A,p}$ is surjective, the group generated by the higher modular points are isomorphic to the $\Q_{p}[\Gn]$-module $V_{(n)}$ or $V_{n}$ depending on the reduction type of $p$ with respect to $E$. In this section, we will prove the following.

\begin{theorem} \thlabel{=}
Let $T_{n}:=\ker\big(\Z_{p}\big[\rfrac{\Gn}{\Bn}\big]\ra \Z_{p}\big)$ be the standard $\Z_{p}[\Gn]$-lattice in $V_{n}$ and let $U$ be a $\Z_{p}[\Gn]$-lattice of $V_{n}$ such that $U^{B_{n}}=T_{n}^{B_{n}}$. Then $U \cong T_{n}$.
\end{theorem}

This will be crucial in the proof of \thref{big} as the saturated group of higher modular points as defined in section $5$ will satisfy the requirements of this theorem.

Let \[V_{(n)}=\Ind_{\Bn}^{\Gn}(\Q_{p})=\Q_{p}\big[\rfrac{\Gn}{\Bn}\big]\] which is a $\Q_{p}[\Gn]$-module of dimension $p^{n+1}+p^{n}$. We can view $V_{(n)}=\text{Maps}(\Pn \ra \Q_{p})$. This decomposes as \[V_{(n)}:=\bigoplus_{i=-1}^{n}W_{i}\] where $W_{i}$ are irreducible $\Q_{p}[\Gn]$-modules as seen in \cite[Theorem 5]{CW09}. These take the form $W_{-1}=\Q_{p},$ $W_{0}=\ker(\Q_{p}\big[\rfrac{G_{0}}{B_{0}}\big]\ra\Q_{p})$ and 
\begin{equation*}
\begin{split}
W_{i}	&=\ker\Big(\Q_{p}\big[\rfrac{G_{i}}{B_{i}}\big]\ra\Q_{p}\big[\rfrac{G_{i-1}}{B_{i-1}}\big]\Big)\\
	&=\bigg\{f:\PR^{1}_{i} \ra \Q_{p}:\sum\limits_{C \supset D}f(C)=0\text{ for all }D \in \PR^{1}_{i-1}\bigg\}
\end{split}
\end{equation*}
for $i \ge 1$. We define the standard $\Z_{p}[\Gn]$-lattice inside $V_{(n)}$ as \[T_{(n)}=\Ind_{\Bn}^{\Gn}(\Z_{p})=\Z_{p}\big[\rfrac{\Gn}{\Bn}\big]=\text{Maps}(\Pn \ra \Z_{p}).\] Define \[V_{n}:=V_{(n)}/W_{-1}=\ker\Big(\Q_{p}\big[\rfrac{\Gn}{\Bn}\big]\ra \Q_{p}\Big)\] and let \[T_{n}:=\ker\Big(\Z_{p}\big[\rfrac{\Gn}{\Bn}\big]\ra \Z_{p}\Big)=\bigg\{f:\Pn \ra \Z_{p}:\sum\limits_{C}f(C)=0\bigg\}\] be the standard $\Z_{p}[\Gn]$-lattice in $V_{n}$. We initially look at the cohomology of $\Bn$ with respect to the lattice $T_{(n)}$.

\begin{proposition} \thlabel{b()}
For all $n \ge 0$, we have $\HC^{1}(\Bn,T_{(n)})=0$.
\end{proposition}

\begin{proof}
By definition, $T_{(n)}=\Ind_{\Bn}^{\Gn}(\Z_{p})$, so we have \[\HC^{1}(\Bn,T_{(n)}) \cong \Ext_{\Z_{p}\Bn}^{1}(\Z_{p},T_{(n)}) \cong \Ext_{\Z_{p}\Gn}^{1}(T_{(n)},T_{(n)}).\] We see by \cite[Corollary 3.3.5 (vi)]{DJB91} that 
\begin{equation*}
\begin{split}
\Ext_{\Z_{p}\Gn}^{1}(T_{(n)},T_{(n)})	&\cong \bigoplus_{\Bn g \Bn}
										\Ext_{\Z_{p}[B_{n} \cap g \Bn g^{-1}]}^{1}(
										\Res^{g \Bn g^{-1}}_{B_{n} \cap 
										g \Bn g^{-1}}(\Z_{p}),\Res^{\Bn}
										_{B_{n} \cap g \Bn g^{-1}}(\Z_{p})) \\
										&\cong \bigoplus_{\Bn g \Bn}
										\Ext_{\Z_{p}[B_{n} \cap g \Bn g^{-1}]}^{1}(
										\Z_{p},\Z_{p})
\end{split}
\end{equation*}
where the sum is taken over the double cosets $\Bn g \Bn$. However, for all such double cosets, we have $\Bn \cap g \Bn g^{-1} \leq \Gn$ and it acts trivially on $\Z_{p}$. Therefore we have 
\begin{equation*}
\begin{split}
\Ext_{\Z_{p}[B_{n} \cap g \Bn g^{-1}]}^{1}(\Z_{p},\Z_{p})	&\cong \HC^{1}(B_{n} \cap 
														g \Bn g^{-1},\Z_{p}) \\
														&\cong \Hom(B_{n} \cap
														 g \Bn g^{-1},\Z_{p})=0,
\end{split}
\end{equation*}
for all double cosets $\Bn g \Bn$. Hence we have $\Ext_{\Z_{p}\Gn}^{1}(T_{(n)},T_{(n)})=0$.
\end{proof}

We then use this for the following.

\begin{corollary} \thlabel{bt}
For all $n \ge 0$, we have $\HC^{1}(\Bn,T_{n})=0$.
\end{corollary}

\begin{proof}
We have the short exact sequence \[0 \ra T_{n} \ra T_{(n)} \xrightarrow{g} \Z_{p} \ra 0\] which induces \[0 \ra T_{n}^{B_{n}} \ra T_{(n)}^{B_{n}} \xrightarrow{g} \Z_{p} \rightarrow \HC^{1}(B_{n},T_{n}) \ra 0\] as $B_{n}$ acts trivially on $\Z_{p}$ and $\HC^{1}(\Bn,T_{(n)})=0$ by \thref{b()}. We know $\Bn$ fixes a point of $\Pn$, say $C$. Then the characteristic function $e_{C}:\Pn \ra \Z_{p}$ in $T_{(n)}^{\Bn}$ is mapped to $1 \in \Z_{p}$ by $g$. Therefore, $g$ is surjective and $\HC^{1}(\Bn,T_{n})=0$.
\end{proof}

We would like to focus on the $\Z_{p}[\Gn]$-lattices $U$ of $V_{n}$ such that $U^{\Bn}=T_{n}^{\Bn}$. All such lattices can be scaled such that they contain $T_{n}$ and are contained in $\frac{1}{p^{k}}T_{n}$ for some $k \ge 1$. We have the exact sequence \[0 \ra T_{n} \ra U \ra \rfrac{U}{T_{n}} \ra 0.\] By \thref{bt}, if $U^{\Bn}=T_{n}^{\Bn}$ then $(\rfrac{U}{T_{n}})^{B_{n}}=0.$ We will initially focus on the $k=1$ case. We now define $\overline{X}:=X \otimes_{\Z_{p}} \F_{p}$ for any $\Z_{p}[\Gn]$-module $X$. We look at the fixed part by $B_{n}$ of the $\F_{p}[\Gn]$-submodules of $\overline{T}_{(n)}:=\F_{p}\big[\rfrac{\Gn}{\Bn}\big]$.

\begin{proposition} \thlabel{not0}
Let $W$ be a non-trivial $\F_{p}[\Gn]$-submodule of $\overline{T}_{(n)}$. Then $W^{B_{n}} \neq 0$.
\end{proposition}

\begin{proof}
Let $W$ be a non-trivial $\F_{p}[\Gn]$-submodule of $\overline{T}_{(n)}$. We will show that the only irreducible submodules of $\overline{T}_{(n)}$ are the trivial module and the Steinberg representation \[W_{st}=\bigg\{f:\Pn \ra \F_{p}: f(C)=t_{D}\text{ for $C \supset D \in \Po$ and $\sum\limits_{D}t_{D}=0$}\bigg\}\] and these modules have non-trivial fixed part by a Borel subgroup.

We will first show these are the only irreducible $\F_{p}[\Gn]$-submodules of $\overline{T}_{(n)}$. Let $W$ be an irreducible $\F_{p}[\Gn]$-module and let \[W_{H}=\rfrac{W}{\langle (h-1)W:\text{ }h \in H\rangle}\] where $H \leq \Gn$. We let $\Bn$ be a Borel subgroup and \[\Pi_{n}=\Bigg\{\begin{pmatrix} a & b \\ 0 & 1 \end{pmatrix} \in \Gn:a \equiv 1 \pmod{p} \Bigg\}.\] This is a Sylow $p$-subgroup of $\Bn$. We will look at the $n=0$ case first. Then $\Pi_{0}$ is generated by $u=(\begin{smallmatrix} 1 & 1 \\ 0 & 1 \end{smallmatrix})$. We see that 
\begin{equation*}
\begin{split}
\Hom_{\F_{p}[\Go]}(W,\overline{T}_{(0)})	&\cong \Hom_{\F_{p}																	[B_{0}]}(W,\F_{p})\\
														&\cong\Hom_{\F_{p}}(W_{B_{0}}
														,\F_{p})\\
														&\cong\Hom_{\F_{p}
														[\rfrac{B_{0}}{\Pi_{0}}]}
														(W_{\Pi_{0}},\F_{p}).
\end{split}
\end{equation*}
We see in \cite[pg. 87]{RG11} that the irreducible $\F_{p}[\Go]$ representations take either the form \[\tau_{j,1}=\Sym^{j-1} \otimes \text{det}^{-\frac{j-1}{2}} \quad \text{or} \quad \tau_{j,2}=\tau_{j,1} \otimes \tau_{1}\] where $\tau_{1}$ is the unique non-trivial $1$-dimensional $\F_{p}[\Go]$-representation and $j \in \{1,3,...,p\}$. We can view $\Sym^{j-1}$ as an $\F_{p}$-subspace of $\F_{p}[x,y]$ of degree $j-1$ homogeneous polynomials where $\GL_{2}(\F_{p})$ acts on $\Sym^{j-1}$ by linear substitution. Then \[(u-1)\Sym^{j-1} \oplus \F_{p}y^{j-1}=\Sym^{j-1}.\] Therefore if we twist by $\det^{-\frac{j-1}{2}}$, we have $\big(\Sym^{j-1} \otimes \text{det}^{-\frac{j-1}{2}}\big)_{\Pi_{0}} \cong \F_{p}y^{j-1} \otimes \text{det}^{-\frac{j-1}{2}}$. This means that if $W$ is an irreducible $\F_{p}[\Go]$-module, then $W$ must have the underlying representation $\tau_{j,1}$ or $\tau_{j,2}$ for $j \in \{1,3,...,p\}$. Hence $W_{\Pi_{0}}$ is a $1$-dimensional irreducible $\F_{p}[\rfrac{\Bo}{\Pi_{0}}]$-module. We see that \[\rfrac{\Bo}{\Pi_{0}} \cong \Bigg\{\begin{pmatrix} a & 0 \\ 0 & 1 \end{pmatrix} \in \Go \Bigg\}\] and $y^{j-1} \in \Sym^{j-1} \otimes \text{det}^{-\frac{j-1}{2}}$ is an eigenvector of $(\begin{smallmatrix} a & 0 \\ 0 & 1 \end{smallmatrix})$ with eigenvalue $a^{-\frac{j-1}{2}}$. We also have \[\tau_{1} \begin{pmatrix} a & 0 \\ 0 & 1 \end{pmatrix}=a^{\frac{p-1}{2}}.\] We see from \cite{RG11} that $\F_{p}$ and $W_{st}$ have underlying representations $\tau_{1,1}$ and $\tau_{p,2}$ respectively. Therefore, $W_{\Pi_{0}}=\F_{p} \iff W=\F_{p}$ or $W=W_{st}$.

For $n \ge 1$, we see that $\Irr_{\F_{p}}(\Gn)=\Irr_{\F_{p}}(\Go)$ where $\Irr_{\F_{p}}(G)$ are the irreducible representations for the group $G$ over $\F_{p}$. As $W$ is an irreducible $\F_{p}[\Gn]$-module, then we can again view $W$ having an underlying representation isomorphic to either $\tau_{j,1}=\Sym^{j-1} \otimes \det^{-\frac{j-1}{2}}$ or $\tau_{j,2}=\tau_{j,1} \otimes \tau_{1}$ for $j \in \{1,3,...,p\}$. We have $\Gn$ acting on it by linear substitution. We also see that \[\Pi_{n} = \Bigg\langle\begin{pmatrix} 1 & 1 \\ 0 & 1 \end{pmatrix}\Bigg\rangle \rtimes \Bigg\{\begin{pmatrix} a & 0 \\ 0 & 1 \end{pmatrix} \in \Gn: a \equiv 1 \pmod{p}\Bigg\}\] which is no longer cyclic. We therefore need to look at the action of the right hand side group on $\Sym^{j-1}$.  Let $\overline{a}$ represent the element $(\begin{smallmatrix} a & 0 \\ 0 & 1 \end{smallmatrix}) \in \Pi_{n}$ where $a \equiv 1 \pmod{p}$. Looking at the $1$-dimensional $\F_{p}[\Gn]$-representation $\det:\Gn \ra \F^{\times}_{p}$, we see that $\det^{-\frac{j-1}{2}}(\overline{a})=1$ as $a \equiv 1 \pmod{p}$. Also, as $\Sym^{j-1}$ can be viewed as degree $j-1$ homogeneous polynomials over $\F_{p}$, then \[\overline{a} \cdot x^{j-1-i}y^{i}=a^{j-1-i}x^{j-1-i}y^{i}=x^{j-1-i}y^{i}\] as $a \equiv 1 \pmod{p}$. Therefore, this shows that $\overline{a}$ fixes $\Sym^{j-1} \otimes \text{det}^{-\frac{j-1}{2}}$ and so \[\big(\Sym^{j-1} \otimes \text{det}^{-\frac{j-1}{2}}\big)_{\Pi_{n}}=\big(\Sym^{j-1} \otimes \text{det}^{-\frac{j-1}{2}}\big)_{H_{n}}\] where $H_{n} \leq \Gn$ is generated by $u=(\begin{smallmatrix} 1 & 1 \\ 0 & 1 \end{smallmatrix})$. Hence we have 
\begin{equation*}
\begin{split}
\Hom_{\F_{p}[\Gn]}(W,\overline{T}_{(n)})	&\cong \Hom_{\F_{p}[B_{n}]}(W,\F_{p})\\
										&\cong\Hom_{\F_{p}}(W_{B_{n}},\F_{p})\\
										&\cong\Hom_{\F_{p}[\rfrac{B_{n}}{\Pi_{n}}]}
										(W_{\Pi_{n}},\F_{p}) \\
										&\cong\Hom_{\F_{p}[\rfrac{B_{n}}{\Pi_{n}}]}
										(W_{H_{n}},\F_{p}). \\
\end{split}
\end{equation*}
We see that if $\overline{a}$ is a representative of a coset in $\rfrac{\Bn}{\Pi_{n}}$, the representation underlying $\Sym^{j-1} \otimes \text{det}^{-\frac{j-1}{2}}$ maps $\overline{a}$ to $a^{-\frac{j-1}{2}}$ and the representation underlying $\big(\Sym^{j-1} \otimes \text{det}^{-\frac{j-1}{2}}\big) \otimes \tau_{1}$ maps $\overline{a}$ to $a^{\frac{p-j}{2}}$. Hence for $W_{H_{n}}=\F_{p}$, we need $j=1$ in the first case or $j=p$ in the second, so $W$ is either $\F_{p}$ or $W_{st}$. 

We finally need to show these irreducible modules have non-trivial fixed part by a Borel subgroup. The trivial $\F_{p}[\Gn]$-module has $\Gn$ acting trivially on it and hence so does every Borel subgroup. Also, for a fixed element of $\Po$, denoted $D_{0}$, then a Borel subgroup acts transitively on the set of $C \not\supset D_{0}$ for $C \in \Pn$. As the number of all elements $C \not\supset D_{0}$ is divisible by $p$, then the fixed part of $W_{st}$ by a Borel subgroup contains the function $f$ such that 
\[    \ f(C)= \left\{
                \begin{array}{ll}
                  1 & \text{if $C \not\supset D_{0}$},\\
                  0 & \text{otherwise},
                \end{array}
              \right.
\]
and hence is non-trivial.
\end{proof}

Therefore, if $U$ is a $\Z_{p}[\Gn]$-lattice inside $\frac{1}{p}T_{n}$ containing $T_{n}$ with $U^{\Bn}=T_{n}^{\Bn}$, then as ${(\rfrac{U}{T_{n}})^{\Bn}=0}$, we can apply \thref{not0} to see that $\rfrac{U}{T_{n}}=0$ meaning $U \cong T_{n}$. 

We now look at the $\Z_{p}[\Gn]$-lattices containing $T_{n}$ that are contained in $\frac{1}{p^{k}}T_{n}$  for some $k \ge 1$. We would like to see when these lattices have non-trivial fixed part by a Borel subgroup.

\begin{proposition} \thlabel{nonzero}
Let $W$ be a non-trivial $\Z/p^{k}\Z[\Gn]$-module contained in $\rfrac{T_{(n)}}{p^{k}T_{(n)}}$ for $k \ge 1$. Then $W^{B_{n}} \neq 0$.
\end{proposition}

\begin{proof}
We do this by induction. The $k=1$ case is \thref{not0} so assume this is true for $k$. We have a short exact sequence \[0 \ra \overline{T}_{(n)} \ra \rfrac{T_{(n)}}{p^{k+1}T_{(n)}} \ra \rfrac{T_{(n)}}{p^{k}T_{(n)}}\ra 0\] so we can view $\overline{T}_{(n)}$ as a $\Z/p^{k+1}\Z[\Gn]$-submodule of $\rfrac{T_{(n)}}{p^{k+1}T_{(n)}}.$ Let $W$ be a non-trivial $\Z/p^{k+1}\Z[\Gn]$-module contained in $\rfrac{T_{(n)}}{p^{k+1}T_{(n)}}$. If $W \cap \overline{T}_{(n)} \neq 0$, then $W$ contains an irreducible $\F_{p}[\Gn]$-module $U$ such that $U^{B_{n}} \neq 0$ by \thref{not0}. Hence we can assume $W \cap \overline{T}_{(n)}=0$. Then by the second isomorphism theorem, we have \[\frac{W+\overline{T}_{(n)}}{\overline{T}_{(n)}}\cong \frac{W}{W \cap \overline{T}_{(n)}} \cong W.\] Hence \[W \cong \frac{W+\overline{T}_{(n)}}{\overline{T}_{(n)}} \leq \frac{\rfrac{T_{(n)}}{p^{k+1}T_{(n)}}}{\overline{T}_{(n)}} \cong \rfrac{T_{(n)}}{p^{k}T_{(n)}}.\] As $W \neq 0$, then by induction, $W^{B_{n}} \neq 0$. 
\end{proof}

We can then finish this section by proving the theorem.

\begin{proof} [Proof of \thref{=}]
We can scale $U$ such that it contains $T_{n}$ and there exists a $k \ge 1$ such that $U \subseteq \frac{1}{p^{k}}T_{n}$. As $U^{\Bn}=T_{n}^{\Bn}$, then we have seen that $(\rfrac{U}{T_{n}})^{\Bn}=0$. Therefore by \thref{nonzero}, $U \cong T_{n}$.
\end{proof} 

\section{DERIVATIVES}

We now look at creating derivative classes associated to the higher modular points with the method proposed by Wuthrich in \cite{CW09} in order to create points of prime power order in certain Selmer groups. To do this, we will use the link created between the higher modular points and representations of $\Gn$ in section 3 as well as \thref{=}. We will show that under the conditions of the prime $p$ with respect to $A$ and $E$ as outlined in \thref{p div}, we understand the type of $\Z_{p}[\Gn]$-lattice the saturated group generated by the higher modular points must be isomorphic to. Further, we can use this to not only show that the order of certain Selmer groups must be divisible by a certain prime power, but also contain an element of prime power order.  

Let $E/\Q$ and $A/\Q$ be elliptic curves of conductors $N_{E}$ and $N_{A}$ respectively. Let $p$ be a prime of good ordinary or multiplicative reduction with respect to $E$ and let \[
    \ F:= \left\{
                \begin{array}{ll}
                  K_{N_{E}} & \text{if $p\nmid N_{E}$},\\
                  K_{\frac{N_{E}}{p}} & \text{if $p || N_{E}$}.
                \end{array}
              \right.
\] We assume that $\tau_{F,A,p}$ is surjective. Hence, if we let $F_{n}$ be the smallest field extension of $F$ such that $\Hn:=\Gal(F_{n}/F)$ acts as scalars on $A[p^{n+1}]$, then we have $\Hn \cong \Gn$.

We define $A_{n}$ to be a non-split Cartan subgroup of $\Gn$. This is a cyclic subgroup of order $p^{n+1}+p^{n}$. We then let $L_{n}$ be the subfield of $F_{n}$ fixed by $A_{n}$. 

\begin{theorem} \thlabel{big}
Let $E/\Q$ and $A/\Q$ be elliptic curves of conductor $N_{E}$ and $N_{A}$ respectively. Let $p>2$ be a prime of multiplicative reduction for $A$ such that $\ord_{p}(\Delta_{A})=1$ and let $F$ be as defined above. Assume that:
	\begin{enumerate}
	\item $A$ does not have potentially good supersingular reduction for any prime of additive reduction;
	\item $E$ has either split multiplicative reduction at $p$ with $p \nmid \ord_{p}(\Delta_{E})$, non-split multiplicative reduction at $p$, or good ordinary non-anomalous reduction at $p$;
	\item The degree of any isogeny of $A$ defined over $\Q$ is coprime to $N_{E}$;
	\item $\tau_{F,A,p}$ is surjective;
	\item $\overline{\rho}_{\Q,E,p}$ is surjective;
	\item Any prime $\ell$ of bad reduction for $E$ and good reduction for $A$ has $a_{\ell}(A)^{2}-4\ell$ square modulo $p$.
	\end{enumerate}
Then there exists an element of order $p^{n}$ in $\Sel^{p^{n}}(E/L_{n})$.
\end{theorem}

In order to show this, we need to look at the splitting of primes in the field extension $F_{n}/L_{n}$. We will prove the following.
\begin{theorem} \thlabel{sc}
Suppose none of the primes of additive reduction for $A$ are potentially good supersingular. Then the extension $F_{n}/L_{n}$ is nowhere ramified. Furthermore, all places above $\infty$, $p$ and $N_{A}$ split completely in the extension as well as all places above a prime $\ell$ dividing $N_{E}$ but not $p$ or $N_{A}$ such that $a_{\ell}(A)^{2}-4\ell$ is a square modulo $p$.
\end{theorem}

We will first need the following lemma.

\begin{lemma}\cite[Lemma 23]{CW09} \thlabel{26}
Let $\nu$ be either a place of ordinary reduction above $p$, an infinite place or a place of potentially multiplicative reduction all with respect to $A$. Then the image of $\overline{\tau}_{F_{\nu},A,p^{n+1}}$ lies in a Borel subgroup of $\Gn$.
\end{lemma}

Therefore, we now need to look at the places of $F$ of bad reduction with respect to $E$ such that they are not above $p$ or $N_{A}$. Let $\nu$ be such a place. Then the Frobenius element $\Fr_{\nu}$ of $\Hn$ generates the decomposition group.

Let $A_{\nu}$ be the reduced elliptic curve at $\nu$ and $R_{\nu}$ the subring of $\End(A_{\nu})$ generated by the Frobenius endomorphism. We define \[u_{\nu}:=\disc(R_{\nu}), \qquad \delta_{\nu}:=0,1\text{ depending on when }u_{\nu}\equiv 0,1 \pmod{4},\] and $b_{\nu}$ the unique positive integer such that $u_{\nu}b_{\nu}^{2}=a_{\nu}(A)^{2}-4q_{\nu}$. We then associate to $\nu$ the integral matrix \[M_{\nu}=\begin{pmatrix} \frac{a_{\nu}(A)+b_{\nu}\delta_{\nu}}{2} & b_{\nu} \\ \frac{b_{\nu}(u_{\nu}-\delta_{\nu})}{4} & \frac{a_{\nu}(A)-b_{\nu}\delta_{\nu}}{2} \end{pmatrix}.\] We see by \cite[Theorem 2.1]{DT02} that as $\nu \nmid p$, then $\nu$ is unramified in $F(A[p^{n+1}])/F$ and $M_{\nu},$ when reduced modulo $p^{n+1}$, represents the conjugacy class of the Frobenius of $\nu$ in $\Gal(F(A[p^{n+1}])/F)$. Thus, we want to see when $M_{\nu}$ reduced modulo $p^{n+1}$ lies in a Borel subgroup. 
\begin{lemma} \thlabel{square}
If $a_{\nu}(A)^{2}-4q_{\nu}$ is a square modulo $p$ then $M_{\nu}$ reduced modulo $p^{n+1}$ lies in a Borel subgroup of $\GL_{2}(\Z/p^{n+1}\Z)$.
\end{lemma}

\begin{proof}
First assume $(b_{\nu},p)=1$. Then $M_{\nu}$ is conjugate to \[\begin{pmatrix} 0 & 1 \\ -q_{\nu} & a_{\nu}(A) \end{pmatrix}\] by \cite[Theorem 2.2]{AOPV09}. This matrix lies in a Borel subgroup if $a_{\nu}(A)^{2}-4q_{\nu}$ is a square modulo $p$.

If $b_{\nu}=p^{j}t$ for some $j \ge 1$ and $(t,p)=1$ then as $\det(M_{\nu})=q_{\nu}$ is not divisible by $p$ then neither is $a_{\nu}(A)$. Therefore, $M_{\nu}$ is conjugate to \[\begin{pmatrix} \frac{a_{\nu}(A)}{2} & p^{j} \\[4pt] \frac{p^{j}t^{2}u_{\nu}}{4} & \frac{a_{\nu}(A)}{2} \end{pmatrix}\] by \cite[Theorem 2.2]{AOPV09}. This matrix lies in a Borel subgroup if $u_{\nu}$ is a square modulo $p$. But as $u_{\nu}b_{\nu}^{2}=a_{\nu}(A)^{2}-4q_{\nu}$, this is true if $a_{\nu}(A)^{2}-4q_{\nu}$ is a square mod $p$.
\end{proof}

Therefore, we need to find conditions for $a_{\nu}(A)^{2}-4q_{\nu}$ to be a square modulo $p$. Let $\ell$ be a prime of good reduction with respect to $A$. Then the arithmetic Frobenius $\Fr_{\ell}$ has characteristic polynomial $x^{2}-a_{\ell}(A)x+\ell$ with roots $\alpha$ and $\beta$. Then we let $t_{n}:=\tr(\Fr_{\ell}^{n})$ for $n \ge 1$ and we have $\ell^{n}=\det(\Fr_{\ell}^{n})$.

\begin{proposition} \thlabel{basecase}
Let $A/\Q$ be an elliptic curve and $p$ a prime. Let $\ell$ be a prime of good reduction with respect to $A$ such that $\ell \neq p$. Then for all $n\ge 1$ \[a_{\ell}(A)-4\ell\text{ is a square mod $p$ if and only if $t_{n}-4\ell^{n}$ is a square mod $p$.}\]
\end{proposition}

We will first look at the cases where $n$ is odd.

\begin{lemma} \thlabel{weird}
For all $n \ge 0$ \[t_{2n+1}^{2}-4\ell ^{2n+1}=(t_{1}^{2}-4\ell)\Bigg((1-2n)\ell^{n}+\sum\limits_{j=1}^{n}\ell^{n-j}t_{j}^{2}\Bigg)^{2}.\]
\end{lemma}

\begin{proof}
Let $\alpha$ and $\beta$ be as before. Then we have \[t_{2n+1}^{2}-4\ell^{2n+1}=(\alpha^{2n+1}-\beta^{2n+1})^{2}.\] Therefore we have 
\begin{equation*}
\begin{split}
\alpha^{2n+1}-\beta^{2n+1}	&=(\alpha-\beta)\Bigg(\sum\limits^{2n}_{j=0}
							\alpha^{j}\beta^{2n-j}\Bigg) \\
							&=(\alpha-\beta)\Bigg(\alpha^{n}\beta^{n}
							+\sum\limits^{n}_{j=1}
							\alpha^{n+j}\beta^{n-j}+\alpha^{n-j}
							\beta^{n+j}\Bigg)\\
							&=(\alpha-\beta)\Bigg((1-2n)\alpha^{n}\beta^{n}
							+\sum\limits^{n}_{j=1}
							(\alpha^{n+j}\beta^{n-j}
							+2\alpha^{n}\beta^{n}+\alpha^{n-j}
							\beta^{n+j})\Bigg)\\
							&=(\alpha-\beta)\Bigg((1-2n)\alpha^{n}\beta^{n}
							+\sum\limits^{n}_{j=1}
							\alpha^{n-j}\beta^{n-j}
							(\alpha^{j}+	\beta^{j})^{2}\Bigg)\\
							&=(\alpha-\beta)\Bigg((1-2n)\ell^{n}
							+\sum\limits^{n}_{j=1}
							\ell^{n-j}t_{j}^{2}\Bigg).\\
\end{split}
\end{equation*}
Squaring both sides gives the equation.
\end{proof}

We then just need to check the case for $n$ even. We first need the following lemma.

\begin{lemma} \thlabel{2n}
For all $n \ge 1$, we have $t_{2n}=t_{n}^{2}-2\ell^{n}$.
\end{lemma}

\begin{proof}
Let $\alpha$ and $\beta$ be as above. Then we have \[t_{2n}=\alpha^{2n}+\beta^{2n}=(\alpha^{n}+\beta^{n})^{2}-2\alpha^{n}\beta^{n}
=t_{n}^{2}-2\ell^{n}. \qedhere \] 
\end{proof}

\begin{proof}[Proof of \thref{basecase}]
We see from \thref{weird} that if $n$ is odd, this is true. If $n$ is even, then we can show by induction and \thref{2n} that if $n=2^{r}m$ with $m$ odd then \[t_{n}^{2}-4\ell^{n}=(t_{m}^{2}-4\ell^{m})\Bigg(\prod_{j=0}^{r-1}t_{2^{j}m}
\Bigg)^{2}. \qedhere\] 
\end{proof}

\begin{proof}[Proof of \thref{sc}]
As $F_{n} \subset F(A[p^{\infty}])$, then it is unramified outside $\infty$, $p$ and $N_{A}$. We see by \thref{26} that the decomposition group of a place in $F$ dividing $\infty \cdot p \cdot N_{A}$ inside $\Hn$ is contained in a Borel. Furthermore, by \thref{square}, the decomposition group of a place $\nu$ in $F$ dividing $N_{E}$ but not $p$ or $N_{A}$ such that $a_{\nu}(A)^{2}-4q_{\nu}$ is a square modulo $p$ inside $\Delta_{n}$ is contained in a Borel. Since any Borel intersects trivially with $A_{n}$ by \cite[Lemma 22]{CW09}, then these places must split completely in $F_{n}/L_{n}$. We then have our result by \thref{basecase}.
\end{proof}

We are now able to prove \thref{big}. 

\begin{proof}[Proof of \thref{big}]
We will follow the proof of \cite[Theorem 21]{CW09}. There is an injection \[\mu:V_{n} \ra E(F_{n}) \otimes \Q_{p}\]\[f \mapsto \sum\limits_{D}f(D) \cdot Q_{A,D}\] where $Q_{A,D}$ are the higher modular points. Let \[S_{n}=\{P \in E(F_{n}):\text{ there exists a $k \ge 0$ such that $p^{k} \cdot P \in \Z_{p}[G_{n}] \cdot Q_{A,D}$}\}\] be the saturated group generated by the higher modular points in $E(F_{n})$. The torsion subgroup of $S_{n}$ is $E(F_{n})[p^{\infty}]$ and so there exists a short exact sequence \[0\ra E(F_{n})[p^{\infty}] \ra S_{n} \ra U_{n} \ra 0\] where $U_{n}$ can be identified as a $G_{n}$-stable lattice in the image of $\mu$ which has no $A_{n}$-fixed elements. If $E(F_{n})[p]=0$, we have $S_{n} \cong U_{n}$ and therefore, by looking at the $p$-adic version of \cite[Proposition 26]{CW09}, $\HC^{1}(A_{n},S_{n})$ has order $p^{n}$. We need to determine when $E(F_{n})[p]=0$.

\begin{proposition} \thlabel{ef0}
If $\overline{\rho}_{\Q,E,p}$ and $\tau_{F,A,p}$ are surjective then for all $n \ge 0 $, $E(F_{n})[p]=0$.
\end{proposition}

\begin{proof}
We first look at the $n=0$ case. As $\overline{\rho}_{\Q,E,p}$ is surjective, then the smallest Galois extension $K$ of $\Q$ such that $E(K)[p] \neq 0$ is $K=\Q(E[p])$. Suppose $E(F_{0})[p] \neq 0$. As $F_{0}/\Q$ is Galois, then we must have $\Q(E[p]) \subseteq F_{0}$. Therefore we have an inclusion of fields $F \subseteq F(E[p]) \subseteq F_{0}$. As $F(E[p])/F$ and $F_{0}/F$ are Galois extensions, then $\Gal(F_{0}/F(E[p])) \trianglelefteq \Gal(F_{0}/F) \cong  \Go$ as $\tau_{F,A,p}$ is surjective. Therefore, $\Gal(F_{0}/F(E[p]))$ must either be trivial, isomorphic to $\Go$ or isomorphic to $\PSL_{2}(\F_{p})$.

If $\Gal(F_{0}/F(E[p]))$ is trivial, then $F(E[p])=F_{0}$ and so ${\Gal(F(E[p])/F) \cong \Go}$. However, this means that $\Go$ is isomorphic to a normal subgroup of $\GL_{2}(\F_{p})$ which cannot happen by \cite[Theorem 4.9]{EA57}, as this would imply a group of that size must be isomorphic to $\SL_{2}(\F_{p})$.

If $\Gal(F_{0}/F(E[p])) \cong \Go$, then $F(E[p])=F$. As $F/\Q$ is a Galois extension, then $\Q(E[p]) \subseteq F$ and so $\mu_{p} \subseteq F$ where $\mu_{p}$ are the $p$-th roots of unity. Therefore, ${\Gal(F(A[p])/F) \leq \SL_{2}(\F_{p})}$ and so ${\Gal(F_{0}/F) \leq \text{PSL}_{2}(\F_{p})}$, which is a contradiction. 

Therefore, $\Gal(F_{0}/F(E[p])) \cong \PSL_{2}(\F_{p})$. Hence, $F(E[p])/F$ must be a quadratic extension such that $\Gal(F(E[p])/F) \cong \im(\chi)$. Here, $\chi$ is the only non-trivial $1$-dimensional representation over $\C$ admitted by $\Go$, which is defined as \[\chi:\Go \xrightarrow{\det} \F_{p}^{\times}/(\F_{p}^{\times})^{2} \ra \{\pm 1\}\] with kernel $\PSL_{2}(\F_{p})$. However, we know $\Gal(F(E[p])/F)$ is isomorphic to a normal subgroup of $\GL_{2}(\F_{p})$ and therefore, by \cite[Theorem 4.9]{EA57}, we must have \[\Gal(F(E[p])/F) \cong \{\pm I_{2}\},\] which is a contradiction. Therefore, $E(F_{0})[p]=0$.

For the $n \ge 1$ case, we let $H_{n}=\Gal(F_{n}/F_{0})$ and $M_{n}=E(F_{n})[p]$. Then we know $H_{n}$ acts on $M_{n}$ and we have already shown $M_{n}^{H_{n}}=0$. However, as $H_{n}$ is a $p$-group and all $H_{n}$-orbits have size greater than $1$ except the orbit representing the trivial element, then $|M_{n}| \equiv 1 \pmod{p}$. Therefore, as $M_{n}$ either has size divisible by $p$ or is trivial, then we must have $M_{n}=0$ and therefore, $E(F_{n})[p]=0$. 
\end{proof}

If we consider the natural inclusion of $S_{n}$ in $E(F_{n}) \otimes \Z_{p}$, then the cokernel $Y_{n}$ is a free $\Z_{p}$-module and we obtain a long exact sequence \begin{equation} 0 \ra E(L_{n}) \otimes \Z_{p} \ra Y_{n}^{A_{n}} \ra \HC^{1}(A_{n},S_{n}) \ra \HC^{1}(A_{n},E(F_{n}))[p^{\infty}]\end{equation} where $Y_{n}^{A_{n}}$ has the same rank as $E(L_{n}) \otimes \Z_{p}$. If we compose the last map with the inflation map we obtain \[\delta_{n}:\HC^{1}(A_{n},S_{n}) \ra \HC^{1}(L_{n},E)[p^{\infty}]\] known as the derivation map. In particular, we have \[\Bigg(\frac{S_{n}}{p^{n}S_{n}}\Bigg)^{A_{n}} \cong \HC^{1}(A_{n},S_{n}).\] We thus call the image of $\delta_{n}$ the \emph{derived classes of higher modular points.} 

\begin{lemma}
The image of $\delta_{n}$ is contained in $\Sh(E/L_{n})[p^{\infty}]$.
\end{lemma}

\begin{proof}
We can follow the method of \cite[Lemma 27]{CW09}. Let $\eta$ be a lift of an element in the image of $\delta_{n}$ in the map \[\HC^{1}(L_{n},E[p^{m}]) \ra \HC^{1}(L_{n},E)[p^{m}]\] for a sufficiently large $m$. The extension $F_{n}/L_{n}$ is unramified at a place $\nu$ outside the set of places in $L_{n}$ above $N_{E}$, $p$ or $\infty$ and so the restriction of $\eta$ to $\HC^{1}(L_{n,\nu},E[p^{m}])$ will lie in $\HC^{1}_{f}(L_{n,\nu},E[p^{m}])$. 

If $\nu$ is a place of $L_{n}$ lying above $\infty$ or $N_{A}$ or is a place that lies above $N_{E}$ but not above $\infty$ or $N_{A}$, then \thref{sc} shows that $\nu$ splits completely in the extension $F_{n}/L_{n}$. Then the restriction of $\eta$ to $\HC^{1}(L_{n,\nu},E)[p^{m}]$ is trivial as it comes from the inflation \[\HC^{1}(F_{n}/L_{n},E(F_{n}))[p^{\infty}] \ra \HC^{1}(L_{n},E)[p^{\infty}].\] Hence, $\eta$ belongs to $\Sel^{p^{m}}(E/L_{n})$.
\end{proof}

We therefore have the map \[\delta_{n}:\HC^{1}(A_{n},S_{n}) \ra \Sh(E/L_{n})[p^{\infty}].\] We now want to determine the source of this mapping. With the conditions that we have on the higher modular points, we have the following.

\begin{proposition}
For all $n \ge 0$, we have $S_{n}^{B_{n}}=T_{n}^{B_{n}}$.
\end{proposition}

\begin{proof}
We show this by induction. For the $n=0$ case, we have $T_{0}^{\Bo} \cong \Z_{p}$, generated by a fixed higher modular point $Q_{A,D_{0}}$. As $Q_{A,D_{0}}$ is not divisible by $p$ in $E(F_{0})$, then $S_{0}^{\Bo}=T_{0}^{\Bo}$. 

Now suppose $n \ge 1$. Let $P \in S_{n}^{\Bn}$ and $k \ge 0$ be minimal such that $p^{k} \cdot P \in T_{n}^{\Bn}$. Assume that $k>0$. We see that $T_{n}^{\Bn}$ is generated by the higher modular points $\{Q_{A,D_{0}},...,Q_{A,D_{n}}\}$ such that $D_{i}$ is a cyclic group of order $p^{i+1}$ and $\Nm(Q_{A,D_{i}})=a_{p}(E) \cdot Q_{A,D_{i-1}}$ for all $1 \leq i \leq n$ as seen in Section \ref{s3}, where $\Nm$ denotes the norm map from $F_{i}^{B_{i}}$ to $F_{i-1}^{B_{i-1}}$. In particular, we see that \[\rfrac{T_{n}^{\Bn}}{T_{n-1}^{B_{n-1}}} = \langle \overline{Q}_{A,D_{n}} \rangle \cong \Z_{p}.\] Therefore, we have $p ^{k} \cdot P=bQ_{A,D_{n}}+R$ for some $b \in \Z_{p}$ and $R \in T_{n-1}^{B_{n-1}}$. We then apply $\Nm$ from $F_{n}^{\Bn}$ to $F_{n-1}^{B_{n-1}}$. We see that $\Nm(Q_{A,D_{n}})=a_{p}(E) \cdot Q_{A,D_{n-1}}$ and $\Nm(R)=pR$. Therefore, we have $p^{k} \cdot \Nm(P)=ba_{p}(E) \cdot Q_{A,D_{n-1}}+pR \in T_{n-1}^{B_{n-1}}$ and $\Nm(P) \in S_{n-1}^{B_{n-1}}=T_{n-1}^{B_{n-1}}$ by induction. As $Q_{A,D_{n-1}}$ is not divisible by $p$ in $E(F_{n-1})$, then $ba_{p}(E)$ is divisible by $p$. However, as the reduction type of $p$ with respect to $E$ is either multiplicative or good ordinary, then $p \nmid a_{p}(E)$ and so $p \mid b$. Therefore, \[R=p\bigg(p^{k-1} \cdot P - \frac{b}{p} \cdot Q_{A,D_{n}}\bigg) \in p \cdot E(F_{n}^{\Bn}).\] We have seen that $E(F_{n})[p]=0$, hence $R$ must already be divisible by $p$ in $E(F_{n-1}^{B_{n-1}})$ and so $R=pR'$ for some $R' \in E(F_{n-1}^{B_{n-1}})$. In particular, we have $R' \in S_{n-1}^{B_{n-1}}=T_{n-1}^{B_{n-1}}$. Therefore, as $p^{k} \cdot P= bQ_{A,D_{n}}+R$, then $p^{k-1} \cdot P=\frac{b}{p} \cdot Q_{A,D_{n}}+R' \in T_{n}^{\Bn}$ which contradicts the minimality of $k$. Therefore, $k=0$ and $P \in T_{n}^{\Bn}$. 
\end{proof}

Using this and \thref{=}, we have $S_{n}\cong T_{n}$. Hence we can look at the cohomology with respect to the standard $\Z_{p}[\Hn]$-lattice in $V_{n}$.

\begin{lemma}
We have $\HC^{1}(A_{n},T_{n}) \cong \Z/p^{n}\Z$.
\end{lemma}

\begin{proof}
Identical to the proof of \cite[Lemma 25]{CW09} which looks at $\Z[\Delta_{n}]$-modules.
\end{proof}

We thus have two options. Firstly, if $\delta_{n}$ is not injective, then from the long exact sequence $(1)$ the rank of $Y_{n}^{A_{n}}$ must be positive and hence so must the rank of $E(L_{n})$. Then $E(L_{n})$ will contribute a copy of $\Z/p^{n}\Z$ in $\Sel^{p^{n}}(E/L_{n})$ as $E(L_{n})[p^{n}]=0$.

If $\delta_{n}$ is injective, then $\Sh(E/L_{n})[p^{\infty}]$ must contain a cyclic subgroup of order $p^{n}$. Therefore, $\Sel^{p^{n}}(E/L_{n})$ must contain an element of order $p^{n}$ by lifting the image of $\delta_{n}$ from $\Sh(E/L_{n})[p^{n}]$.
\end{proof}

We can look further at \thref{big} by taking semistable $A/\Q$ and assuming $\overline{\rho}_{\Q,A,p}$ is surjective. We will require the following.

\begin{proposition} 
Let $A/\Q$ be semistable and $p>3$ be a prime such that $\overline{\rho}_{\Q,A,p}$ is surjective. Then $\tau_{F,A,p}$ is surjective.
\end{proposition}

\begin{proof}
If $\overline{\rho}_{F,A,p}$ is surjective, then as $A$ is semistable, $\rho_{F,A,p}$ must be surjective and hence $\tau_{F,A,p}$ is also. Therefore, we need to show that for all primes $p>3$, $\overline{\rho}_{F,A,p}$ is surjective. 

We initially look at the case when $p$ is a prime of split multiplicative reduction with respect to $A$. We have seen by construction that $F \subseteq \Q(A[M])$ where 
\[
    \ M= \left\{
                \begin{array}{ll}
                  N_{E} & \text{if $p\nmid N_{E}$},\\
                  \frac{N_{E}}{p} & \text{if $p || N_{E}$},
                \end{array}
              \right.
\]  
so $(M,p)=1$. Let $K=\Q(A[p]) \cap \Q(A[M])$ and $\nu$ be a place of $K$ above $p$. Then 
\begin{align*}
K_{\nu}	&\cong \Q_{p}(A[p]) \cap \Q_{p}(A[M]) \\
		&=\Q_{p}(\mu[p],\sqrt[p]{q_{A}}) \cap \Q_{p}(\mu[M],\sqrt[M]{q_{A}})
\end{align*}
as $p$ is a prime of split multiplicative reduction. Here, $\mu[n]$ denotes the group of $n^{\operatorname{th}}$ roots of unity and $q_{A} \in \Q_{p}^{\times}$ is a $p$-adic analytic uniformiser at $\infty$ such that $\ord_{p}(q_{A})=-\ord_{p}(j(A))>0$. We want to show $K_{\nu}=\Q_{p}$. As ${(M,p)=1}$, then $\Q_{p}(\mu[p]) \cap \Q_{p}(\mu[M])=\Q_{p}$. Therefore, \[\Gal(\Q_{p}(\mu[pM])/\Q_{p}) \cong G \leq (\Z/M\Z)^{\times} \times (\Z/p\Z)^{\times}.\] Let $L=\Q_{p}(\mu[pM])$. Then we see by Kummer theory that $\Gal(L(\sqrt[p]{q_{A}})/L)$ is a subgroup of $\Z/p\Z$ and $\Gal(L(\sqrt[M]{q_{A}})/L)$ is a subgroup of $\Z/M\Z$. As $(M,p)=1$, then we must have $L(\sqrt[p]{q_{A}}) \cap L(\sqrt[M]{q_{A}})=L$ and so the Galois group of ${L(\sqrt[pM]{q_{A}})=L(\sqrt[p]{q_{A}},\sqrt[M]{q_{A}})}$ over $L$ is the direct product of these Galois groups. We also see that \[L \cap \Q_{p}(\mu[M],\sqrt[M]{q_{A}})=\Q_{p}(\mu[M])\] and \[L \cap \Q_{p}(\mu[p],\sqrt[p]{q_{A}})=\Q_{p}(\mu[p]).\] Therefore, looking at the Galois groups, we have \[\Gal(L(\sqrt[pM]{q_{A}})/\Q_{p}(\mu[p],\sqrt[p]{q_{A}})) \cong \Gal(\Q_{p}(\mu[M],\sqrt[M]{q_{A}})/\Q_{p}).\] As $L(\sqrt[pM]{q_{A}})$ is the compositum of $\Q_{p}(\mu[M],\sqrt[M]{q_{A}})$ and $\Q_{p}(\mu[p],\sqrt[p]{q_{A}})$, then this implies \[\Q_{p}(\mu[p],\sqrt[p]{q_{A}}) \cap \Q_{p}(\mu[M],\sqrt[M]{q_{A}})=\Q_{p}.\] Therefore, $K_{\nu}=\Q_{p}$. Now, let $H=\Gal(\Q(A[p])/K)$. Then as $K_{\nu}=\Q_{p}$, the inertia group, $I_{\nu}$ of $H$ must contain the set of matrices of the form $(\begin{smallmatrix} * & 0 \\ 0 & 1 \end{smallmatrix})$ under the image of $\overline{\rho}_{\Q,A,p}$ by \cite[pg. 277, Corollaire]{JPS72}. As $K/\Q$ is Galois, then $H \trianglelefteq \GL_{2}(\F_{p})$ and so as $p \ge 5$, then either $H \supseteq \SL_{2}(\F_{p})$ or $H \subseteq Z_{p}$ by \cite[Theorem 4.9]{EA57}. Therefore, as $I_{\nu}$ contains the set of matrices of the form $(\begin{smallmatrix} * & 0 \\ 0 & 1 \end{smallmatrix})$, then $H \supseteq \SL_{2}(\F_{p})$. But then the mapping $\det:H \ra \F_{p}^{\times}$ is surjective and has kernel $\SL_{2}(\F_{p})$ meaning $H=\GL_{2}(\F_{p})$. Therefore, we have $\Q(A[p]) \cap \Q(A[M])=\Q$. However, as $F \subseteq \Q(A[M])$, then $\Q(A[p]) \cap F=\Q$ and as $\overline{\rho}_{\Q,A,p}$ is surjective, then so is $\overline{\rho}_{F,A,p}$.

Now suppose $p$ is a prime of non-split multiplicative reduction for $A$. Let $K=\Q(A[p]) \cap \Q(A[M])$ and let $\nu$ be a place of $K$ above $p$. As $p$ is a prime of non-split multiplicative reduction for $A$, there exists an unramified quadratic extension $L'/\Q_{p}$ such that $A$ acquires split multiplicative reduction over $L'$. In particular, we can follow the same procedure as in the split case with $L'$ to show that $L'(A[p]) \cap  L'(A[M])=L'$. This now ensures that $K_{\nu} \in \{\Q_{p},L'\}$ and so $K_{\nu}/\Q_{p}$ is an unramified extension. Therefore, if $H=\Gal(\Q(A[p])/K)$, then the image of the inertia subgroup, $I_{\nu}$ of $H$ under $\overline{\rho}_{\Q,A,p}$ must contain the set of matrices of the form $(\begin{smallmatrix} * & 0 \\ 0 & 1 \end{smallmatrix})$ by \cite[pg. 277, Corollaire]{JPS72}. Therefore, we again see that as $H \trianglelefteq \GL_{2}(\F_{p})$, then $H=\GL_{2}(\F_{p})$ and so $\Q(A[p]) \cap \Q(A[M])=\Q$. Therefore, identically to the split case, this means that $\overline{\rho}_{F,A,p}$ is surjective and this completes the proof.
\end{proof}

Therefore, when $A/\Q$ is semistable, we can reduce \thref{big} to the following.

\begin{theorem} \thlabel{semibig}
Let $E/\Q$ and $A/\Q$ be elliptic curves of conductor $N_{E}$ and $N_{A}$ respectively. Let $p>3$ be a prime of multiplicative reduction for $A$ such that $\ord_{p}(\Delta_{A})=1$ and let 
\[
    \ F_{n}:= \left\{
                \begin{array}{ll}
                  K_{p^{n+1}N_{E}} & \text{if $p\nmid N_{E}$},\\
                  K_{p^{n}N_{E}} & \text{if $p || N_{E}$}.
                \end{array}
              \right.
\] for $n \ge -1$ with $F:=F_{-1}$. Assume that:
	\begin{enumerate}
	\item $A$ is semistable;
	\item $E$ has either split multiplicative reduction at $p$ with $p \nmid \ord_{p}(\Delta_{E})$, non-split multiplicative reduction at $p,$ or good ordinary non-anomalous reduction at $p$;
	\item The degree of any isogeny of $A$ defined over $\Q$ is coprime to $N_{E}$;
	\item $\overline{\rho}_{\Q,A,p}$ is surjective;
	\item $\overline{\rho}_{\Q,E,p}$ is surjective;
	\item Any prime $\ell$ of bad reduction for $E$ and good reduction for $A$ has $a_{\ell}(A)^{2}-4\ell$ square modulo $p$.
	\end{enumerate}
Then there exists an element of order $p^{n}$ in $\Sel^{p^{n}}(E/L_{n})$.
\end{theorem}

\section{EXAMPLES}

We now look at applying \thref{semibig} to a few examples. Due to the size of the field extensions $L_{n}/\Q$, it would be very difficult to verify the calculations in the examples. For information on specific elliptic curves, we obtained our data from \cite{lmfdb}.

\begin{example}
Let \[E:y^{2}+y=x^{3}-x^{2}-x-2\] be the elliptic curve with Cremona label $143a1$ and \[A:y^{2}+y=x^{3}+x^{2}-x\] be the elliptic curve with Cremona label $35a3$. Then if we let $p=7$, this is a prime of good ordinary non-anomalous reduction for $E$. We also see that as $\Delta_{A}=-35$, then $7$ is a prime of multiplicative reduction for $A$ such that $\ord_{7}(\Delta_{A})=1$. The only $\Q$-isogenies on $A$ are of degree $3$ or $9$, which are coprime to $N_{E}$ and both $\overline{\rho}_{\Q,A,7}$ and $\overline{\rho}_{\Q,E,7}$ are surjective. Also,
\begin{equation*}
\begin{split}
a_{11}(A)^{2}-4 \cdot 11=-35 \equiv 0 \pmod{7}\\
a_{13}(A)^{2}-4 \cdot 13=-27 \equiv 1 \pmod{7}
\end{split}
\end{equation*}
which are both squares modulo $7$. Therefore by \thref{semibig}, there exists an element of order $7^{n}$ in $\Sel^{7^{n}}(E/L_{n})$ for all $n \ge 1$.
\end{example}
We also see that $E$ doesn't necessarily have to be a semistable elliptic curve.
\begin{example}
Let \[E:y^{2}+xy=x^{3}-x^{2}-5\] be the elliptic curve with Cremona label $45a1$ and \[A:y^{2}+xy+y=x^{3}+x^{2}\] be the elliptic curve with Cremona label $15a8$. Then if we let $p=5$, then this is a prime of non-split multiplicative reduction for $E$ and a prime of multiplicative reduction for $A$ such that $\ord_{5}(\Delta_{A})=1$. The only $\Q$-isogenies on $A$ are of degree $2$, $4$, $8$ or $16$, which are coprime to $N_{E}$ and both $\overline{\rho}_{\Q,A,5}$ and $\overline{\rho}_{\Q,E,5}$ are surjective. Therefore by \thref{semibig}, there exists an element of order $5^{n}$ in $\Sel^{5^{n}}(E/L_{n})$ for all $n \ge 1$.
\end{example}
We can even vary the curves $E$ and $A$ if they both have prime conductor. However, we first require the following.
\begin{proposition} \thlabel{atleast1}
Let $E/\Q$ be an elliptic curve of prime conductor $p$. Then there exists an elliptic curve, $A/\Q$, in its isogeny class such that $\ord_{p}(\Delta_{A})=1$.
\end{proposition}
\begin{proof}
Let $E_{0}$ be the strong Weil curve in the isogeny class of $E$. Then by \cite[Proposition 2]{DW09}, there are three cases. Either
\begin{itemize}
\item $\#E_{0}(\Q)=1$,
\item $\#E_{0}(\Q)=2$,
\item $\#E_{0}(\Q)>2$.
\end{itemize}
For the first case, we have $E=E_{0}$ and as $p$ is a prime of multiplicative reduction for $E$, then ${\ord_{p}(\Delta_{E})=c_{p}=1}$ by \cite[Proposition 2]{DW09}, where $c_{p}$ is the Tamagawa number of $E$.  For the second case, $p$ is of the form $u^{2}+64$ for some integer ${u \equiv 3 \pmod{4}}$ by \cite[Proposition 2]{DW09}. The elliptic curve $E$ is one of the following two isogenous curves
\begin{equation*}
\begin{split}
&E_{0}:y^{2}+xy=x^{3}-\frac{u+1}{4}x^{2}+4x-u \\
&E_{1}:y^{2}+xy=x^{3}-\frac{u+1}{4}x^{2}-x
\end{split}
\end{equation*}
where $\Delta_{E_{1}}=p$. Finally, if $\#E_{0}(\Q)>2$, then the conductor of $E_{0}$ will be either $11, 17, 19$ or $37$. But then $E_{0}$ and hence $E$ will be isogenous to one of the curves with Cremona label $11a2$, $17a3$, $19a3$ or $37b2$, which all have $p$ strictly dividing their discriminant.
\end{proof}
\begin{example}
Let $E/\Q$ be an elliptic curve of prime conductor $p$ and let $A/\Q$ be an elliptic curve of conductor $p$ such that $\ord_{p}(\Delta_{A})=1$, which exists by \thref{atleast1}. We see in \cite[Proposition 2]{DW09} that if $p$ is a prime of split multiplicative reduction for $E$ then $p \nmid \ord_{p}(\Delta_{E})$ and $\overline{\rho}_{\Q,E,p}$ is surjective for all $p$. We also have that $\overline{\rho}_{\Q,A,p}$ is surjective as $p \ge 11$. In particular, $A$ has no $p$-isogeny defined over $\Q$. Hence, by \thref{semibig}, $\Sel^{p^{n}}(E/L_{n})$ contains an element of order $p^{n}$ for all $p$ prime and $n \ge 1$.
\end{example}
We also show a specific example for $E$ with conductor $q$ when there exists an $A$ of conductor $pq$ for $p$, $q$ prime such that $p \neq q$.
\begin{example}
Let $a,b \in \Z$ such that at least one of $a$ or $b$ is not divisible by $3$ and \[a^{12}-9a^{8}b+27a^{4}b^{2}-27b^{3}\] is not a square. Then by \cite[Theorem 3.10]{HJ12} there exists infinitely many $m \in \Z$ such that \[A_{m}:y^{2}+y=x^{3}+ax^{2}+bx+m\] is an elliptic curve where $|\Delta_{A_{m}}|=pq.$ Here, $p$ is an odd prime and $q$ is either $1$ or an odd prime different from $p$. Pick $m$ such that $p,q \ge 11$ and there exists an elliptic curve $E/\Q$ of conductor $q$ such that $p$ is a prime of good ordinary non-anomalous reduction. Then as $A_{m}$ is semistable, $\overline{\rho}_{\Q,A_{m},p}$ and $\overline{\rho}_{\Q,A_{m},q}$ are surjective, which ensures that there doesn't exist a $\Q$-isogeny of $A_{m}$ of degree $q$. In particular, $p$ is a prime of multiplicative reduction such that $\ord_{p}(\Delta_{A_{m}})=1$. Also, $\overline{\rho}_{\Q,E,p}$ is surjective as $p \ge 11$.  Therefore, by \thref{semibig}, $\Sel^{p^{n}}(E/L_{n})$ contains an element of order $p^{n}$ for all $n \ge 1$.
\end{example}

\subsection*{Acknowledgments}
I would like to thank Chris Wuthrich for all his help on this topic.

This work was supported by the Engineering and Physical Sciences Research Council (Grant No. EP/N50970X/1).

\bibliography{Kolyvagin_Derivatives_of_Modular_Points_on_Elliptic_Curves}
\bibliographystyle{abbrv}

\end{document}